\documentclass{amsart}

\usepackage[latin1]{inputenc}
\usepackage[T1]{fontenc}	

\usepackage{amsmath}
\usepackage{amssymb}
\usepackage{graphicx}
\usepackage{url}
\usepackage{mathrsfs}
\usepackage{stmaryrd}
\usepackage{listings}

\newtheorem {theo} {Theorem}
\newtheorem {prop} {Proposition}
\newtheorem {defi} {Definition}
\newtheorem {algo} {Algorithm}

\def\Z{\mathbb{Z}}

\DeclareMathAlphabet{\nim}{U}{bbold}{m}{n}

\title{Nimbers are inevitable}
\author{Julien Lemoine - Simon Viennot}

\begin{document}

\renewcommand\labelitemi{$\ast$}

\begin{abstract}
This article concerns the resolution of impartial combinatorial games, and in particular games that can be split in sums of independent positions. We prove that in order to compute the outcome of a sum of independent positions, it is always more efficient to compute separately the nimbers of each independent position than to develop directly the game tree of the sum. The concept of nimber is therefore inevitable to solve impartial games, even when we only try to determinate the winning or losing outcome of a starting position. We also describe algorithms to use nimbers efficiently and finally, we give a review of the results obtained on two impartial games: Sprouts and Cram.
\end{abstract}

\maketitle


\section{Introduction}

The games considered in this article are \emph{combinatorial games}: two players play alternately, with perfect knowledge of the current state of the game, and there is no room for chance. We will restrain our discussion to \emph{impartial} combinatorial games: from a given position, the same moves are available to both players.

Moreover, the theorems and algorithms of this article apply only to impartial combinatorial games in their \emph{normal} version, where the first player who cannot move loses, and not to the \emph{misère} version, where the first player who cannot move wins\footnote{The analysis of impartial games in misère version is much more complicated, notably because the algorithms described in this article cannot be applied.}.

We will focus our attention particularly on \emph{splittable} impartial games, in which some of the positions can be split in sum of independent positions. Our purpose is to \emph{solve} these games, i.e. to find which player has a winning strategy and to compute it explicitly. In section 2, we review some background notions on impartial games, illustrating them with the games of Sprouts and Cram, and in section 3, we give some insight on the central concept of nimber.

In section 4, we develop the main result of this article: we prove that nimbers are necessary when we try to compute the outcome of a splittable impartial game. In section 5, we detail algorithms to use nimbers efficiently and finally, in section 6, we present the results obtained on the games of Sprouts and Cram.

\section{Background}

\subsection{Sprouts and Cram}

The algorithms described in this paper can be applied to any impartial combinatorial game played in the normal version, and we have chosen two well-known games for our computations: Sprouts and Cram.

The game of Sprouts starts with a given number of spots drawn on a sheet of paper. Alternately, the players draw a line between two spots (possibly the same spot), and add a spot anywhere on this line. The lines cannot cross each other, and a given spot cannot be used in more than 3 lines.

\begin{figure}[ht]
\centering
\includegraphics[scale=0.5]{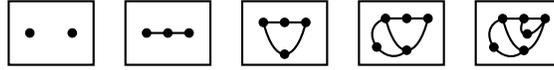}
\caption{Example of a Sprouts game, starting with 2 spots (the second player wins).}
\end{figure}

The first article about Sprouts is from Martin Gardner \cite{mg67} in 1967. A detailed presentation can be found in Winning Ways \cite{ww01}.

The game of Cram is played on a board with very simple rules: players alternately fill two adjacent cells with a domino, until one of them cannot play anymore. A description and some interesting analysis can also be found in \cite{ww01}. 

\begin{figure}[ht]
\centering
\includegraphics[scale=0.5]{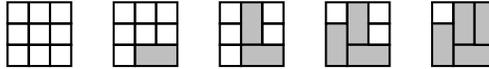}
\caption{Example of a Cram game on a $3\times 3$ board (the second player wins).}
\end{figure}

\subsection{Splittable positions}

Sprouts and Cram are \emph{splittable} games, because some of the positions can be split in a sum of independent positions. When a player moves in such a position, the move can only affect one of the component of the sum, leaving the others untouched.

\begin{figure}[ht]
\centering
\includegraphics[scale=0.5]{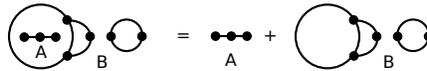}
\caption{Splittable Sprouts position.}
\label{decoupable_sprouts}
\end{figure}

For example, the position of Sprouts on figure \ref{decoupable_sprouts} is splittable. The spots at the interface between regions A and B cannot be used anymore, and any further move must be done inside the region A (without affecting B) or inside the region B (without affecting A).

\begin{figure}[ht]
\centering
\includegraphics[scale=0.5]{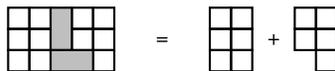}
\caption{Splittable Cram position.}
\label{decoupable_cram}
\end{figure}

Figure \ref{decoupable_cram} gives another example. The position is obtained after playing two moves in a Cram game on a $3\times 5$ board. The position is splittable, because the two components are independent.

\subsection{Game tree}

We call \textit{game tree of a position $\mathscr{P}$} the tree where nodes are the positions obtained by playing moves in $\mathscr{P}$, and in which two positions $\mathscr{P}_1$ and $\mathscr{P}_2$ are linked by an edge if $\mathscr{P}_2$ is an \emph{option} of $\mathscr{P}_1$ (i.e. when $\mathscr{P}_2$ can be reached from $\mathscr{P}_1$ in one move).

\begin{figure}[ht]
\centering
\includegraphics[scale=0.33]{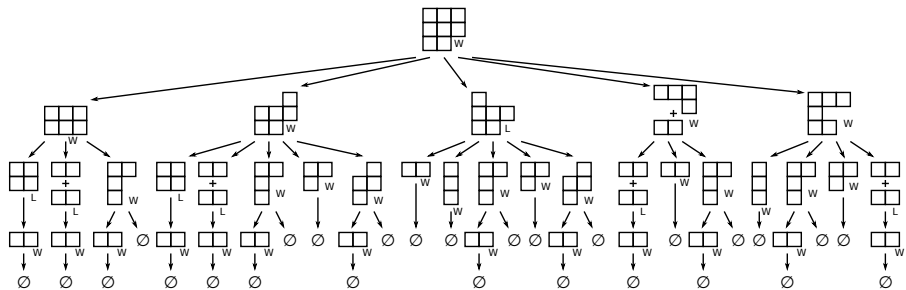}
\caption{Game tree of a Cram position.}
\label{arbre_cram}
\end{figure}

The game tree of figure \ref{arbre_cram} has been obtained by identifying similar positions respectively to symmetry, and deleting isolated cells (since they cannot be used in any further move).

\subsection{Outcome of a position}
\label{issuePosition}

The \emph{outcome} of a position is ``W'' (Win) if, from this position, the player to move has a winning strategy. Otherwise, the outcome is ``L'' (Loss). Positions whose outcome is ``W'' are said to be \emph{winning} and those whose outcome is ``L'' are said to be \emph{losing}.

It is possible to determine recursively the outcome of a position from its game tree. If a position $\mathscr{P}$ has an option which is losing, then $\mathscr{P}$ is winning. Otherwise, all options of $\mathscr{P}$ are winning, and $\mathscr{P}$ is losing. Finally, since the player who cannot move loses, the leaves (terminal positions) are losing.

The outcome of all the positions of figure \ref{arbre_cram} have been indicated.

\subsection{Solution tree}
\label{arbreSolution}
The definition of the outcome of a position shows that it is sufficient to find only one losing option in order to prove that a node is winning. It implies that it is possible to determine the outcome of the root of the tree (the starting position) without knowing the outcome of all the positions of the tree.

\begin{figure}[ht]
\centering
\includegraphics[scale=0.5]{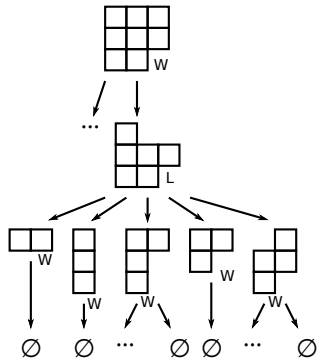}
\caption{Solution tree for a Cram position.}
\label{arbre_cram_mini}
\end{figure}

On figure \ref{arbre_cram_mini}, we have selected a subset of the nodes from figure \ref{arbre_cram}, which is sufficient to prove that the root is losing. There are 3 winning nodes for which it was not necessary to compute the outcome of all the options.

Such a subset of the complete game tree will be called a \emph{solution tree} for the root. A solution tree is a graphic representation of what is also called a ``winning strategy''. If the root is winning, like on figure \ref{arbre_cram_mini}, then the first player has a winning strategy. If the root is losing, the winning strategy is on the contrary for the second player.

The main goal of this article is to describe efficient methods to \emph{solve} splittable impartial games, i.e. to compute a solution tree for the starting positions of these games.

\subsection{Outcome of a sum of positions}

When a position is split in a sum of independent components, it is sometimes possible to compute the outcome of components separately and deduce the result of the sum by using the following result (which can be proved easily by induction):

\begin{prop}
\label{ThSomme}
The sum of two losing positions is losing. The sum of a losing position and a winning position is winning.
\end{prop}

This result allows us to reduce the size of the solution tree when a splittable position is computed. It was used for example by Applegate, Jacobson and Sleator in 1991 to compute the outcome of Sprouts positions \cite{ajs91}. However, it should be noted that it indicates nothing if both components are winning. To determine the outcome of the sum with separate computations even in this case, we need to use the concept of \emph{nimber}.

\section{Nimber}
\label{Nimber}

\subsection{The game of Nim}
\label{JeuNim}

The game of Nim is played with heaps of objects, for example matches. A move consists in removing some of the matches from a single heap, and when the game is played in the normal version, the player that removes the last match wins (because the other player cannot play anymore).

A Nim-heap with n matches will be denoted $\nim{n}$. The position $\nim{7}+\nim{5}+\nim{4}+\nim{2}$ is then composed of 4 heaps with 7, 5, 4 and 2 matches (respectively). For example, the player to move could choose to remove 3 matches in the second heap, which would lead to the position $\nim{7}+\nim{2}+\nim{4}+\nim{2}$. Or he could remove all the matches of the third heap, obtaining $\nim{7}+\nim{5}+\nim{0}+\nim{2}$.

Since a move is restricted to a single heap, the heaps are independent components and the game of Nim is a splittable game. A position from the game of Nim is then the sum of its heaps, which each are an independent component.

The resolution of Nim has been first described by Bouton, in 1902 \cite{bo02}. The method uses the $\oplus$ operator (\emph{bitwise exclusive or}), which we will call \emph{Nim-sum}. To compute the Nim-sum of two integers, we can write them in a binary form, and add the bits with the addition of $\Z/2\Z$ ($0+0=0$, $0+1=1$ and $1+1=0$). For example, $9\oplus 12$ can be written in binary form $1001\oplus 1100$, which gives $0101$ (binary form). Back to the decimal usual notation, we obtain: $9\oplus 12=5$.

The solution of Bouton can be stated as follows:

\begin{theo}[Bouton]
A sum of Nim-heaps has the same outcome as the Nim-sum of the heaps.
\end{theo}

Since the outcome of a single heap is L when it is empty, and W if there are still some matches left (just take all the matches), the losing positions of the complete game of Nim are those for which the Nim-sum of the heaps is 0.

Going back to our first example, it means that the position $\nim{7}+\nim{5}+\nim{4}+\nim{2}$ is winning, because $7\oplus 5\oplus 4\oplus 2=4$. The winning moves are those that leave your opponent in a losing situation, and you should then remove matches so that you get $\nim{3}+\nim{5}+\nim{4}+\nim{2}$ (since $3\oplus 5\oplus 4\oplus 2=0$), or $\nim{7}+\nim{1}+\nim{4}+\nim{2}$, or $\nim{7}+\nim{5}+\nim{0}+\nim{2}$.

\subsection{Indistinguishability}

We will say that two positions $\mathscr{P}_1$ and $\mathscr{P}_2$ are indistinguishable, and will denote $\mathscr{P}_1\sim\mathscr{P}_2$ if, for any position $\mathscr{P}$, the sums of positions $\mathscr{P}_1+\mathscr{P}$ and $\mathscr{P}_2+\mathscr{P}$ have the same outcome. In this definition, the positions $\mathscr{P}_1$, $\mathscr{P}_2$ and $\mathscr{P}$ can be taken from any impartial combinatorial game.

\begin{figure}[ht]
\centering
\includegraphics[scale=0.5]{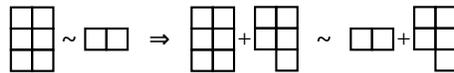}
\caption{A Cram position simplified using indistinguishability.}
\label{cram_equiv}
\end{figure}

This theoretical concept could be useful for practical computations. If a complicated position is known to be indistinguishable from a simpler one, we could replace the complicated one by the simple one in any sum appearing in the computation. For example, the figure \ref{cram_equiv} shows how we could accelerate the computation of figure \ref{decoupable_cram}, if we already knew that the position on the left is indistinguishable of a simpler position (with only two cells). This basic idea will be used in a more sophisticated way with the nimbers.

\subsection{Nimber}

Here is the main result of the theory of impartial games:

\begin{theo}[Sprague-Grundy]
Any position of an impartial game played in the normal version is indistinguishable of some Nim-heap, called its  \emph{nimber}.
\end{theo}

A proof can be found, for example, in Winning Ways \cite{ww01}. The following propositions can be deduced immediately:

\begin{prop}
Let $\mathscr{P}$ a position.
\begin{itemize}
 \item $\mathscr{P}$ is losing $\Leftrightarrow$ the nimber of $\mathscr{P}$ is $\nim{0}$.
 \item $\mathscr{P}$ is winning $\Leftrightarrow$ the nimber of $\mathscr{P}$ is $\geq\nim{1}$.
\end{itemize}
 \end{prop}

\begin{prop}
\label{equiv_couple}
Let $\mathscr{P}$ a position.
\begin{itemize}
 \item $\mathscr{P}\sim\nim{n}\Leftrightarrow \mathscr{P}+\nim{n}$ is losing.
 \item $\mathscr{P}\nsim\nim{n}\Leftrightarrow \mathscr{P}+\nim{n}$ is winning.
\end{itemize}
 \end{prop}

The nimber\footnote{The nimber is also called \emph{number} or \emph{function of Sprague-Grundy}.} has a practical interest in the case of a sum of independent positions. Indeed, it is possible to compute the nimber of the sum from the nimbers of the components, with the Nim-sum.

For example, on figure \ref{decoupable_cram}, the nimber of the first component is $\nim{1}$, and the nimber of the other is $\nim{0}$, so the nimber of the sum is $\nim{1}$ (since $1\oplus 0=1$), and we can then deduce than the sum is winning. Note that this result could have been deduced with proposition \ref{ThSomme}. Now, let us consider figure \ref{decoupable_sprouts}. The nimber of each component is $\nim{2}$, so the nimber of the sum is $\nim{0}$ (since $2\oplus 2=0$), and we can deduce that the sum is losing, which was not possible with only proposition \ref{ThSomme}.

Noting that $m\oplus n=0 \Leftrightarrow m=n$, we obtain:
\begin{prop}
\label{p1p2nim}
Let $\mathscr{P}_1$ and $\mathscr{P}_2$ two positions.
\begin{itemize}
 \item $\mathscr{P}_1+\mathscr{P}_2$ is losing $\Leftrightarrow$ the nimbers of $\mathscr{P}_1$ and $\mathscr{P}_2$ are equal.
 \item $\mathscr{P}_1+\mathscr{P}_2$ is winning $\Leftrightarrow$ the nimbers of $\mathscr{P}_1$ and $\mathscr{P}_2$ are different.
\end{itemize}
 \end{prop}

To complete this review of the concept of nimber, we need to explain how to compute the nimber of a position that is not a sum of independent components. We will use the following definition:

\begin{defi}
We define the $Mex$ (minimum excluded value) of a set of integers as the least positive integer that is not included in the set.
\end{defi}

For example, by applying this definition to nimbers, we obtain: $Mex(\nim{1}, \nim{4}) = \nim{0}$, $Mex(\nim{0}, \nim{1}, \nim{2}, \nim{5}) = \nim{3}$.

The following proposition allows us to compute recursively the nimber of a position, knowing that the nimber of a terminal position is always $\nim{0}$:

\begin{prop}
The nimber of a position is equal to the $Mex$ of the nimbers of its options.
\end{prop}

\begin{figure}[ht]
\centering
\includegraphics[scale=0.33]{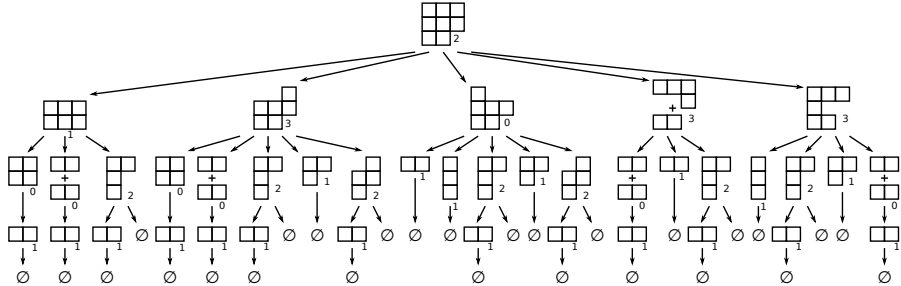}
\caption{Nimbers of the game tree of a Cram position.}
\label{arbre_cram_nim}
\end{figure}

We can now determine the nimbers of all the positions of figure \ref{arbre_cram}, as shown in the figure \ref{arbre_cram_nim}.

Since the definition of $Mex$ uses all the options of a given position, it is sometimes assumed that it is needed to develop the complete game tree of a position in order to compute its nimber. For this reason, nimbers are sometimes considered to require too much running time when we are only trying to compute the winning or losing outcome of a given starting position. The next section will show that, surprisingly, this is not the case at all.

\section{Nimbers are inevitable}

\subsection{Elementary computation of the outcome of a sum of positions}
\label{elem}
The most simple way to compute the outcome of a sum of positions $\mathscr{P}_1+\mathscr{P}_2$ is to consider the complete sum $\mathscr{P}_1+\mathscr{P}_2$ as a single position, and compute the outcome of the options.

\begin{figure}[ht]
\centering
\includegraphics[scale=0.5]{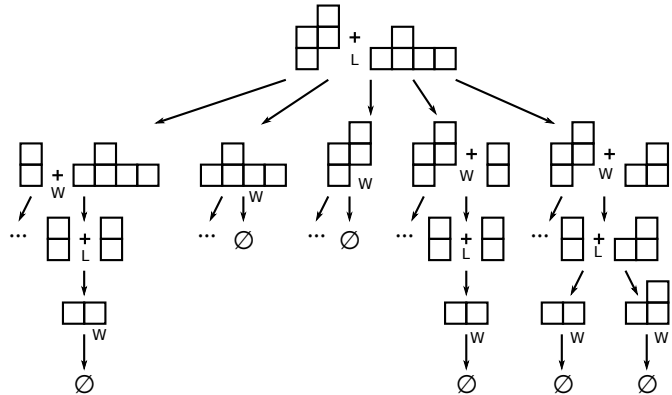}
\caption{Solution tree for a sum of independent positions.}
\label{arbre_cram_somme}
\end{figure}

The figure \ref{arbre_cram_somme} is an example of this method. Note that all the sums appearing in this solution tree -- in particular the starting position -- are of the kind W+W, and proposition \ref{ThSomme} is then of no help to simplify the computation of the outcome.

\subsection{Inevitability of the nimbers}

We can now state the main theoretical result of this article, which shows that even when we are only computing the outcome of a sum of positions, it is always more efficient to compute separately the nimbers of the components.

\begin{theo}[Inevitability of the nimbers]
Let us suppose that we have computed a solution tree for the sum $\mathscr{P}_1+\mathscr{P}_2$ of two independent positions. Then, without computing any other node, we can determine the nimber of one component, and whether the nimber of the other component is equal or different.
\end{theo}

\begin{proof}
This result is proved by induction:
\begin{itemize}
\item Terminal case: if $\mathscr{P}_1=\varnothing$ and $\mathscr{P}_2=\varnothing$, $\mathscr{P}_1+\mathscr{P}_2$ is losing, and the nimber of both $\mathscr{P}_1$ and $\mathscr{P}_2$ is $\nim{0}$.
\item Induction: \begin{itemize}
\item Case 1: if $\mathscr{P}_1+\mathscr{P}_2$ is winning, then one of the option is losing, for example of the form $\mathscr{P}_1^i+\mathscr{P}_2$. In that case, the nimber of $\mathscr{P}_2$ is known by induction hypothesis. And since $\mathscr{P}_1+\mathscr{P}_2$ is winning, the nimbers of $\mathscr{P}_1$ and $\mathscr{P}_2$ are different (with proposition \ref{p1p2nim}).
\item Case 2: if $\mathscr{P}_1+\mathscr{P}_2$ is losing, then all the options are winning. In particular, all the options of the form $\mathscr{P}_1^i+\mathscr{P}_2$ are winning. Either we know the nimber of $\mathscr{P}_2$, or we know the nimbers of all $\mathscr{P}_1^i$ by induction hypothesis, from which we can deduce the nimber of $\mathscr{P}_1$. And since $\mathscr{P}_1+\mathscr{P}_2$ is losing, we conclude that the nimbers of $\mathscr{P}_1$ and $\mathscr{P}_2$ are equal (with proposition \ref{p1p2nim} again).
\end{itemize}
\end{itemize}
\end{proof}

\begin{figure}[ht]
\centering
\includegraphics[scale=0.5]{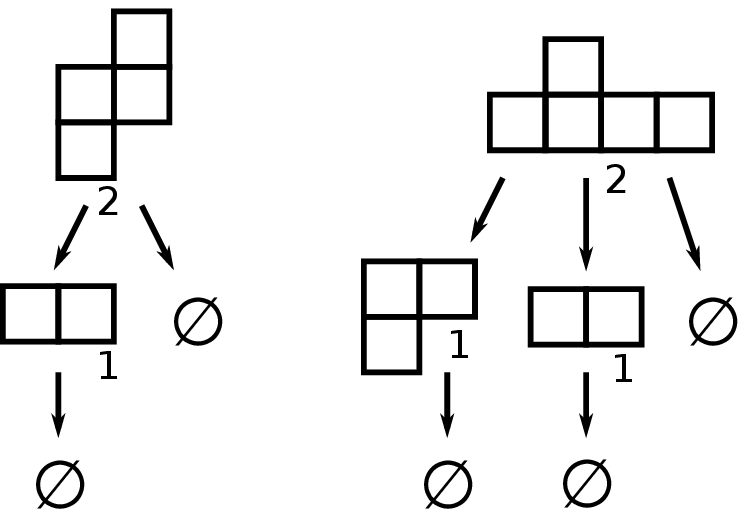}
\caption{Computation of the nimbers of the sum components.}
\label{arbre_cram_somme_nim}
\end{figure}

This result shows that computing the nimbers of the components of the sum (more precisely computing the nimber of one component, and whether the other component has the same nimber or not) does not require to compute more positions than the elementary computation described in paragraph \ref{elem}. And on the contrary, in most cases, it requires strictly less nodes. For example, instead of developing the tree of figure \ref{arbre_cram_somme}, it is more efficient to compute separately the nimbers of each independent component as in the figure \ref{arbre_cram_somme_nim}. The nimber of both positions is $\nim{2}$, so the nimber of the sum is $\nim{0}$ (since $2\oplus 2=0$), from which we deduce finally that the sum is losing.

Despite the simplicity of the demonstration, this result is surprising. Up to now, it was widely assumed that the inherent complexity of the concept of nimber made it unsuitable when trying to compute only the outcome of a sum of very complicated positions\footnote{See for example the discussion at the end of the article of Applegate et al. \cite{ajs91}.}. Yet, the theorem of inevitability proves that the elementary computation of the sum contains at least the computation of one nimber ! It is then always more efficient to use nimbers as an intermediate step in order to compute the outcome of a sum of positions, and that, whatever complicated the positions of the sum are.

\subsection{Computation of the nimber from a partial knowledge of the game tree} 
\begin{figure}[ht]
\centering
\includegraphics[scale=0.5]{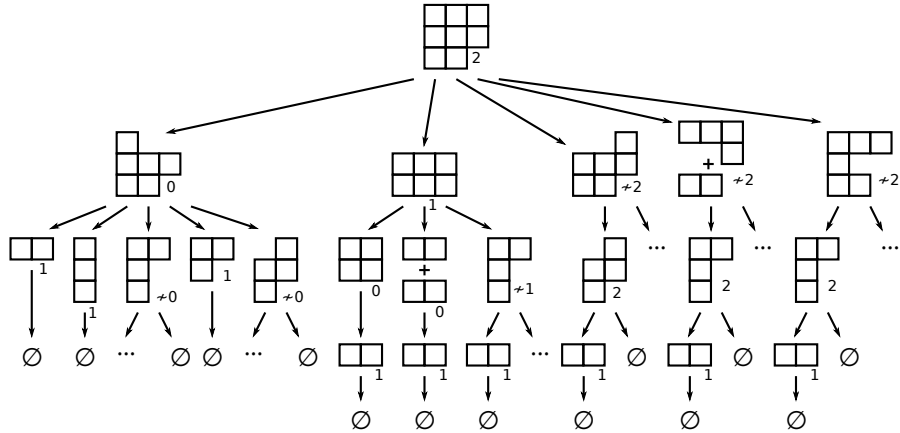}
\caption{Part of the game tree sufficient to compute the nimber of the root.}
\label{arbre_cram_partiel_nim}
\end{figure}

The theorem of inevitability shows that the computation of the nimber of one component is inevitable. However, it does not imply that all the ways of computing this nimber are equivalent: a direct use of the rule of the $Mex$, for example, would cause the development of the complete game tree, which is extremely costly in running time in case of a complicated position. In fact, it is possible, with the algorithms described in the next section, to compute the nimber of a position without developing the whole game tree. We give in figure \ref{arbre_cram_partiel_nim} a sub-tree of figure \ref{arbre_cram_nim}, which is sufficient to determine the nimber of the root. Note in particular that for some nodes, we only prove that their nimber is different from a given value.

\section{Computation algorithms}
\label{algo_calcul}

\subsection{Reformulating nimber computations as outcome computations}

In this section, we will detail algorithms that use the nimber to accelerate the computation of the outcome of sums of positions. We have to deal with two different kinds of computations: computations of outcomes, and computations of nimbers. But the proposition \ref{equiv_couple} shows that it is equivalent to compute that $\mathscr{P}\sim\nim{n}$ (i.e. to compute that the nimber of $\mathscr{P}$ is $\nim{n}$), or to compute that the outcome of $\mathscr{P}+\nim{n}$ is L. Similarly, it is equivalent to compute that $\mathscr{P}\nsim\nim{n}$, or that the outcome of $\mathscr{P}+\nim{n}$ is W. In this way, we can compute whether the nimber of a position is $\nim{n}$ or not, simply by computing the outcome of $\mathscr{P}+\nim{n}$.

In the concrete implementation, we will represent $\mathscr{P}+\nim{n}$ by a \emph{couple} $(\mathscr{P},\nim{n})$. We call $\mathscr{P}$ the \emph{position part} of the couple, and $\nim{n}$ the \emph{nimber part}.

The options of a couple $(\mathscr{P},\nim{n})$ are of two kinds:
\begin{itemize}
\item those of the position part, of the form $(\mathscr{P}^i,\nim{n})$ where $\mathscr{P}^i$ is an option of $\mathscr{P}$.
\item those of the nimber part, of the form $(\mathscr{P},\nim{i})$ with $\nim{i}<\nim{n}$.
\end{itemize}

The computation starts on the couple $(\mathscr{P},\nim{0})$. Indeed, this couple has no option in the nimber part, so we can identify the options of $(\mathscr{P},\nim{0})$ and $\mathscr{P}$. Note that this is true for any position $\mathscr{P}$ and not only the starting position: we can identify $(\mathscr{P},\nim{0})$ and $\mathscr{P}$, and in particular they have the same outcome. 

\subsection{Computation of the outcome of a couple}

The key-point in the case of a splittable game is to check whether $\mathscr{P}$ is splittable or not before computing recursively the outcome of $(\mathscr{P},\nim{n})$. If the position is splittable, we use algorithm \ref{algo_split} described thereafter.

\begin{algo}[Recursive computation of the outcome of a couple]
\label{algo_nimberPart}
~\\To compute the outcome of the couple $(\mathscr{P},\nim{n})$:
\begin{itemize}
\item  If $\mathscr{P}$ is splittable, compute the outcome of the couple with algorithm \ref{algo_split}, otherwise:
\item  For each option $(\mathscr{P}^i,\nim{n})$ of the position part and each option  $(\mathscr{P},\nim{i})$ of the nimber part,
compute the outcome of the option with algorithm \ref{algo_nimberPart}.\\If the option is losing , return ``W''.
\item If all the options are winning, return ``L''.
\end{itemize}
\end{algo}

Let us note that in the case of a non-splittable position $\mathscr{P}$, and by using $\nim{n}=\nim{0}$ as the nimber part, the algorithm is the same as the classical algorithm used to compute the outcome of $\mathscr{P}$ .

\subsection{Computation of the outcome of a sum}

To compute the outcome of a couple when the position part is splittable in a sum of the form $(\mathscr{P}_1+...+\mathscr{P}_k,\nim{n})$, we first compute the nimbers of all the components except one, with algorithm \ref{algo_012} described thereafter. Then, we merge the nimbers with the Nim-sum. The couple is therefore reduced to a couple of the form $(\mathscr{P}_k, \nim{n}'$), without any sum in the position part, and we compute its outcome with algorithm \ref{algo_nimberPart}.

\begin{algo}[Computation of the outcome of a sum]
\label{algo_split}
~\\To compute the outcome of a couple $(\mathscr{P}_1+...+\mathscr{P}_{k},\nim{n})$:
\begin{itemize}
\item For $j$ from $1$ to $k-1$, compute $\nim{n}_j$, the nimber of $\mathscr{P}_j$, with algorithm \ref{algo_012}.
\item Compute $\nim{n}'=\nim{n}_1+...+\nim{n}_{k-1}+\nim{n}$ with the Nim-sum.
\item Compute the outcome of $(\mathscr{P}_k, \nim{n}'$) with algorithm \ref{algo_nimberPart}, and return the obtained value.
\end{itemize}
\end{algo}

\subsection{Computation of the nimber of the position}

Lastly, we need to explain how to compute the nimber of a position, which was necessary in the previous algorithm. The principle is simply to try the nimbers in increasing order: $\nim{0}$, then $\nim{1}$, then $\nim{2}$,... until we find the correct value.

\begin{algo}[Computation of the nimber of a position]
\label{algo_012}
~\\To compute the nimber of a $\mathscr{P}$:
\begin{itemize}
\item Initialise $\nim{n}$ to $\nim{0}$.
\item While the computation of the outcome of $(\mathscr{P},\nim{n})$ with algorithm \ref{algo_nimberPart} returns ``W'', increment $\nim{n}$.
\item Return the final value of $\nim{n}$.
\end{itemize}
\end{algo}

The returned value is the nimber of $\mathscr{P}$, since the loop ends when $(\mathscr{P},\nim{n})$ is found losing.

It should be noted that this algorithm has the advantage of avoiding useless computations, because the computation of $\mathscr{P}\sim\nim{n}$ contains the computation of $\mathscr{P}\nsim\nim{k}$ for $\nim{k}<\nim{n}$. This comes from the fact that if we have proved that $\mathscr{P}\sim\nim{n}$, i.e. we have proved that the couple $(\mathscr{P},\nim{n})$ is losing, then we have proved that each option of this couple is winning. In particular, for any $\nim{k}<\nim{n}$, the option $(\mathscr{P},\nim{k})$ is winning, which means that $\mathscr{P}\nsim\nim{k}$.

\subsection{Game tree traversal}

The efficiency of the algorithms described above depends on the path that we choose in the game tree. In the case of the classical algorithm for computing the outcome of a position, it is sufficient to find a losing option in order to prove that a position is winning. Therefore, the choice of the option that we compute first is important: if we choose a winning option before a losing one, there will be useless computations. And if there is more than one losing option, it is better to choose the ``easiest'' one first, in order to obtain a smaller solution tree, and to obtain it faster.

Such a choice is also needed in the algorithms of the previous section. Firstly, in algorithm \ref{algo_nimberPart}: just as in the classical computation, if the couple $(\mathscr{P},\nim{n})$ is winning, it is better to search the game tree from the easiest losing option first. But a choice also occurs in algorithm \ref{algo_split}: if the couple $(\mathscr{P}_1+...+\mathscr{P}_{k},\nim{n})$ is winning, the nimber of one component will not be computed (in the description of the algorithm, it is $\mathscr{P}_{k}$, but we can choose any of the components). The choice of this component will affect the speed of the computation.

Of course, these choices are not easy to perform, because we do not know which option is the losing one - if we knew it there would be no need for computation ! It should be noted that even if the nimber plays a central role in the algorithms, using the notion of couple enables us to keep some kind of similarity with the classical algorithm to compute the outcome. It follows that most of the usual methods (Depth-First, or Best-First, like the PN-search) used to search game trees in an efficient order can be used in combination with the algorithms of this paper, with some adaptations.

\section{Results}

\subsection{Game of Sprouts}

We have applied the algorithms described in the previous section to the game of Sprouts, which allowed us to compute the outcome of the game up to 32 starting spots and some sparse values up to 47 spots. The following table indicates, for a given number of starting spots $p$, the number of losing couples\footnote{We store only losing couples, i.e. positions whose nimber is known, in order to save memory.} stored in the solution tree obtained at the end of the computation (after pruning useless positions). All the computed outcomes support the ``Sprouts conjecture'': the position with $p$ starting spots is losing if and only if $p=0, 1$ or $2$ modulo 6.

It is interesting to note that before the introduction of the algorithms described in this article, the biggest known outcome was $p=11$, with more than 100,000 losing positions at the end of the computation \cite{ajs91}, whereas we are now able to compute the same position with only 140 losing couples. The algorithms of this article are particularly efficient in the case of Sprouts because splittable positions are extremely frequent, and appear even in the upper part of the game tree.

\begin{figure}[ht]
\centering
\begin{minipage}{0.13\textwidth}
\begin{tabular}{ |c|c| }
\hline
$p$ & size\tabularnewline
\hline
0 & 0\tabularnewline
1 & 1\tabularnewline
2 & 3\tabularnewline
3 & 6\tabularnewline
4 & 16\tabularnewline
5 & 38\tabularnewline
6 & 64\tabularnewline
\hline
\end{tabular}
\end{minipage}
\begin{minipage}{0.15\textwidth}
\begin{tabular}{ |c|c| }
\hline
$p$ & size\tabularnewline
\hline
7 & 103\tabularnewline
8 & 205\tabularnewline
9 & 63\tabularnewline
10 & 140\tabularnewline
11 & 140\tabularnewline
12 & 475\tabularnewline
13 & 577\tabularnewline
\hline
\end{tabular}
\end{minipage}
\begin{minipage}{0.15\textwidth}
\begin{tabular}{ |c|c| }
\hline
$p$ & size\tabularnewline
\hline
14 & 1580\tabularnewline
15 & 3252\tabularnewline
16 & 1068\tabularnewline
17 & 471\tabularnewline
18 & 3233\tabularnewline
19 & 3630\tabularnewline
20 & 4051\tabularnewline
\hline
\end{tabular}
\end{minipage}
\begin{minipage}{0.16\textwidth}
\begin{tabular}{ |c|c| }
\hline
$p$ & size\tabularnewline
\hline
21 & 9270\tabularnewline
22 & 5706\tabularnewline
23 & 2837\tabularnewline
24 & 9316\tabularnewline
25 & 9229\tabularnewline
26 & 18567\tabularnewline
27 & 59117\tabularnewline
\hline
\end{tabular}
\end{minipage}
\begin{minipage}{0.16\textwidth}
\begin{tabular}{ |c|c| }
\hline
$p$ & size\tabularnewline
\hline
28 & 14813\tabularnewline
29 & 3414\tabularnewline
30 & 58363\tabularnewline
31 & 58365\tabularnewline
32 & 58204\tabularnewline
33 & ?\tabularnewline
34 & 21107\tabularnewline
\hline
\end{tabular}
\end{minipage}
\begin{minipage}{0.16\textwidth}
\begin{tabular}{ |c|c| }
\hline
$p$ & size\tabularnewline
\hline
35 & 18812\tabularnewline
... & ?\tabularnewline
40 & 45782\tabularnewline
41 & 48890\tabularnewline
... & ?\tabularnewline
47 & 54542\tabularnewline
... & ?\tabularnewline
\hline
\end{tabular}
\end{minipage}
\caption{Results obtained on the game of Sprouts.}
\end{figure}

The details about the optimisations specific to the game of Sprouts can be found in our article from 2007 \cite{lv07} and the update published in 2010 \cite{lv10}.

\subsection{Game of Cram}

We have also applied the same algorithms to the game of Cram, and we present here the results obtained up to this point.

There exists a symmetry strategy on boards of even$\times$even dimensions, which are losing (and hence their nimber is $\nim{0}$), and similarly on boards of even$\times$odd dimensions, which are winning. But this strategy does not apply to odd$\times$odd boards, and tells nothing about the nimber of even$\times$odd boards (except that this nimber is not $\nim{0}$). In the tables below, we have indicated in parenthesis the results that does not need any computation because of the symmetry strategy. Moreover, we have indicated with ``$-$'' the $n\times m$ board where $n>m$, because the value is the same as on $m\times n$ boards, with a simple symmetry argument.

As far as we know, the best results known so far were those of Martin Schneider in 2009 \cite{ms09}, which we have indicated in the tables with a ``*''.

\begin{figure}[ht]
\centering
\begin{tabular}{|*{16}{c|}}
\hline
3 &  4 &  5 &  6 &  7 &  8 &  9 & 10 & 11 & 12 & 13 & 14 & 15 & 16 & 17 & 18\\ \hline
*$\nim{0}$ & *$\nim{1}$ & *$\nim{1}$ & *$\nim{4}$ & *$\nim{1}$ & *$\nim{3}$ & *$\nim{1}$ & $\nim{2}$ & $\nim{0}$ &  $\nim{1}$ & $\nim{2}$ & $\nim{3}$ &  $\nim{1}$ & $\nim{4}$ & $\nim{0}$ & $\nim{1}$\\ \hline
\end{tabular}
\caption{Results obtained on $3\times n$ boards.}
\end{figure}

\begin{figure}[ht]
\centering
\begin{tabular}{|*{7}{c|}}
\hline
  &  4 &  5 &  6 &  7 &  8 &  9 \\ \hline
4 &($\nim{0}$) & *$\nim{2}$ &($\nim{0}$) &  $\nim{3}$ & ($\nim{0}$)&  $\nim{1}$\\ \hline
5 &  $-$ & *$\nim{0}$ &  $\nim{2}$ & *$\nim{1}$ &  $\nim{1}$ &  W \\ \hline
6 &  $-$ &  $-$ & ($\nim{0}$)& >$\nim{3}$& ($\nim{0}$)& (W)\\ \hline
7 &  $-$ &  $-$ &  $-$ &  W & (W)&    \\ \hline
\end{tabular}
\caption{Results obtained on $n\times m$ boards, with $n\geq 4$ and $m\geq 4$.}
\end{figure}

The new results on boards with both dimensions greater than 4 are the nimbers of the 4$\times$7, 4$\times$9, 5$\times$6, and 5$\times$8 boards, and the winning outcome of the 5$\times$9 and 7$\times$7 boards. With algorithm \ref{algo_012}, we have also been able to compute that the nimber of the 6$\times$7 board is strictly greater than 3, but without achieving to compute the exact value up to now.

In the case of the game of Cram, we have noted experimentally that a slight increase in the board dimensions creates a huge increase of the computation difficulty. This is due of course to the exponential increase of the number of possible board positions, but also to the fact that the greater the board dimensions, the later the splitting of the positions occur in the game tree. The optimisations specific to the game of Cram will be the subject of a future article.

\section*{Conclusion}

Since the discovery of the Sprague-Grundy theorem, the nimbers have been used successfully to analyse a number of impartial games, in particular the numerous variants of Nim, like the octal games. However, in the case of very intricate games, like Sprouts or Cram, the nimbers were usually considered to consume too much running time, and were not or almost not used to compute the outcome of the starting positions.

The theorem presented in this article shows on the contrary that the use of nimbers in impartial splittable games is inevitable, even when we are only trying to compute the outcome of a starting position, because the elementary computation of the outcome of a sum of positions indeed computes the nimber of one of the component. Algorithms using nimbers efficiently have been applied successfully to Sprouts and Cram, two impartial splittable games, and nimbers play a central part in the results obtained.

\bibliography{nimber_en}
\bibliographystyle{amsplain}

\end{document}